\documentclass[]{article}
\begin{document}
\newtheorem{proposition}{Proposition}[section]
\newtheorem{definition}{Definition}[section]
\newtheorem{lemma}{Lemma}[section]

\title{\bf Invariant Algebras}
\author{Keqin Liu\\Department of Mathematics\\The University of British Columbia\\Vancouver, BC\\
Canada, V6T 1Z2}
\date{January, 2011}
\maketitle

\begin{abstract} We introduce invariant algebras and representation$^{(c_1, \dots, c_8)}$ of algebras, and give many ways of constructing Lie algebras, Jordan algebras, Leibniz algebras, 
pre-Lie algebras and left-symmetric algebras in an invariant algebras.
\end{abstract}

\medskip
In this paper, we introduce invariant algebras and representation$^{(c_1, \dots, c_8)}$ of algebras. The main property of invariant algebras is that  an invariant algebra carries $14$ associative products, which can be used to  construct many different algebras inside an invariant algebra. We use this paper to give $6$ ways of constructing a Lie algebra structure, $6$ ways of constructing a Jordan algebra structure, $4$ ways of constructing a Leibniz algebra structure, $14$ ways of constructing a pre-Lie algebra structure and $14$ ways of constructing a left-symmetric algebra structure inside an invariant algebra. Based on the main property of invariant algebras, we introduce representation$^{(c_1, \dots, c_8)}$ of algebras in the last section of this paper. Using representation$^{(c_1, \dots, c_8)}$ of algebras is a natural way of developing representation theory for some non-associative algebras like  Leibniz algebras, pre-Lie algebras and left-symmetric algebras.

\medskip
Throughout, all vector spaces are vector spaces over a field $\mathbf{k}$.

\medskip
\section{Basic Definitions}

\medskip
The notion of invariant algebras introduced in this paper is given in the following

\begin{definition}\label{def1.1} Let $A$ be an associative algebra with an idempotent $q$. The {\bf (right) invariant algebra)} $(A,\, q)$ induced by the idempotent $q$ is defined by
\begin{equation}\label{eq01.1}
(A,\, q):=\{\, x, \,|\, \mbox{$x\in A$ and $qxq=qx$}\,\}.
\end{equation}
\end{definition}

The most important example of invariant algebras is the linear invariant algebra over a vector space. Let $V$ be a vector space, let $End(V)$ be the associative algebra of all linear transformations from $V$ to $V$, and let $I$ (or $I_V$) be the identity linear transformation of $V$. If $W$ is a subspace of $V$, then an linear transformation $q$ satisfying
\begin{equation}\label{eq01.2}
q(W)=0\quad\mbox{and}\mbox (q-I)(V)\subseteq W
\end{equation}
is an idempotent because
$$
q^2(v)-q(v)=q(q-I)(v)\in q(W)=0\quad\mbox{for all $v\in V$.}
$$
An element $q$ of $End(V)$ satisfying (\ref{eq01.2}) is called a {\bf $W$-idempotent}. 

\medskip
Note that 
\begin{equation}\label{eq01.3}
(End(V), q)=\{\, f \,|\, \mbox{$f\in End(V)$ and $f(W)\subseteq W$}\,\},
\end{equation}
where $W$ is a subspace of a vector space $V$ and $q$ is a $W$-idempotent. The invariant algebra $(End(V), q)$ is called the {\bf linear invariant algebra over $V$} induced by the $W$-idempotent $q$, which consists of all linear transformations of $V$ having $W$ as an invariant subspace by (\ref{eq01.3}).

\medskip
Let $(A, q)$ be an invariant algebra. The set
\begin{equation}\label{eq01.7}
(A, q)^{ann}:= \{\, qx-x \,|\, x\in (A, q)\,\}
\end{equation}
is clearly a subspace of $(A, q)$. For $a$, $x\in (A, q)$, we have
$$(qx-x)a=q(xa)-(xa)\in (A, q)^{ann}$$
and
$$a(qx-x)=q(-aqx+ax)-(-aqx+ax)\in (A, q)^{ann}.$$
Hence, the subspace $(A, q)^{ann}$ is an ideal of $(A,q)$, which is called the {\bf annihilator} of $(A, q)$. For $x\in (A,q)$, we have
$$xq-x=q(-xq+x)-(-xq+x)\in (A, q)^{ann},$$
which implies that
\begin{equation}\label{eq01.8}
(A, q)^{ann}\supseteq \{\, xq-x\,|\, x\in (A,q)\,\}\bigcup\{\, qx-xq \,|\, x\in (A,q)\,\}.
\end{equation}

\begin{definition}\label{def01.1} Let $(A,q_{_A})$ and $(B,q_{_B})$ be invariant algebras. A map $\phi: (A,q_{_A}) \to (B,q_{_B})$ is called an {\bf invariant homomorphism} if
$$
\phi(x+y)=\phi(x)+\phi(y), \quad \phi(xy)=\phi(x)\phi(y)\quad\mbox{for $x$, $y\in (A,q_{_A})$}
$$
and
$$ \phi(1_{_A})=1_{_B}, \qquad \phi(q_{_A})=q_{_B},$$
\end{definition}
where $1_{_A}$ and $1_{_B}$ are the identities of $A$ and $B$, respectively. A bijective invariant homomorphism is called an {\bf invariant isomorphism}.

\medskip
The next proposition shows that any invariant algebra can be embedded in a linear invariant algebra.

\begin{proposition}\label{pr01.1} If $(A,q)$ is an invariant algebra, then the map
$$\phi: a\mapsto a_{_L} \quad\mbox{for $a\in (A,q)$}$$
is an injective invariant homomorphism from $(A,q)$ to the linear invariant algebra 
$(End(A,q), q_{_L})$ over $(A,q)$ induced by the $(A, q)^{ann}$-idempotent $q_{_L}$, where $a_{_L}$ is the left multiplication defined by
$$ a_{_L}(x):=ax \quad\mbox{for $x\in (A.q)$.}$$
\end{proposition}

\medskip
\section{Hu-Liu Products}

\medskip
We now define Hu-Liu products in the following 

\begin{definition}\label{def01.2} Let $(A, q)$ be an invariant algebra over a field $\mathbf{k}$. For $1\le i\le 14$, the {\bf $i$-th Hu-Liu product} $\circ_i$ is defined by
\begin{eqnarray}
\label{eq01.13}x\circ_{_1} y  :&=& kyx+hyqx,\\
\label{eq01.14}x\circ_{_2} y  :&=& kyx-kyqx+qxy,\\
\label{eq01.15}x\circ_{_3}y  :&=& yx+kyqx-yxq,\\
\label{eq01.16}x\circ_{_4}y  :&=& yx+kqxy-yxq,\\
\label{eq01.17}x\circ_{_5}y  :&=& kyx-kyqx+qyx,\\
\label{eq01.18}x\circ_{_6}y  :&=& yx+kqyx-yxq,\\
\label{eq01.19}x\circ_{_7}y  :&=& kyx+xqy-kyxq,\\
\label{eq01.20}x\circ_{_8}y  :&=& kxy+hxqy,\\
\label{eq01.21}x\circ_{_9}y  :&=& xy-xqy+kqxy,\\
\label{eq01.22}x\circ_{_{10}}y  :&=& xy-xqy+kqyx,\\
\label{eq01.23}x\circ_{_{11}}y  :&=& xy-xyq+kqxy,\\
\label{eq01.24}x\circ_{_{12}}y  :&=& xy+kyqx-xyq,\\
\label{eq01.25}x\circ_{_{13}}y  :&=& xy-xyq+kqyx,\\
\label{eq01.26}x\circ_{_{14}}y  :&=& xy+kxqy-xyq,
\end{eqnarray} 
where $k$, $h$ are fixed scalars in the field $\mathbf{k}$ and $x$, $y\in (A, q)$.
\end{definition}

\medskip
For $i=1$ or $8$, $\circ_{_i}$ is also denoted by $\circ_{_{i,k, h}}$ to indicate that $\circ_{_1}$ and $\circ_{_8}$ depend on the scalars $k$ and $h$. Similarly, for $7\ge i\ge 2$ or 
$14\ge i\ge 9$, $\circ_{_i}$ is also denoted by $\circ_{_{i,k}}$.

\medskip
The next proposition gives the main property of  $14$ Hu-Liu products.

\begin{proposition}\label{pr01.2} If $(A, q)$ is an invariant algebra over a field $\mathbf{k}$, then the $i$-th Hu-Liu product $\circ_{_i}$ satisfies the associative law:
\begin{equation}\label{eq01.27}
(x\circ_{_i}y)\circ_{_i}z=x\circ_{_i}(y\circ_{_i}z),
\end{equation}
where $1\le i\le 14$ and $x$, $y$, $z\in (A, q)$.
\end{proposition}

\medskip
\section{Lie Algebras}

If $(A, q)$ is an invariant algebra, then $(A, q)$ can be made into a Lie algebra by the well-known square bracket $[\, , \,]$, where $[\, , \,]$ is defined by
\begin{equation}\label{eq01.73}
[x, y]:=xy-yx \quad\mbox{for $x$, $y\in (A, q)$.}
\end{equation}
Except the ordinary square bracket (\ref{eq01.73}), there are other six square brackets which  also make an invariant algebra into a Lie algebra. We now introduce the six square brackets in the following

\begin{definition}\label{def01.4} Let $(A,q)$ be an invariant algebra over a field $\mathbf{k}$. For $1\le i\le 7$, the {\bf $i$-th square bracket} $[\, ,\, ] _i$ is defined by
\begin{eqnarray}
\label{eq01.74} [x, y ]_1 :&=& qxy-qyx,\\
\label{eq01.75} [x, y ]_2 :&=& xqy-yqx,\\
\label{eq01.76} [x, y ]_{3} :&=& xy-yx+kxqy-kyqx,\\
\label{eq01.77} [x, y ]_{4} :&=& xy-yx-xqy+yqx-kqxy+kqyx,\\
\label{eq01.78} [x, y ]_{5} :&=& xy-yx-xyq+yxq-kqxy+kqyx,\\
\label{eq01.79} [x, y ]_{6} :&=& xy-yx-xyq+yxq+kxqy-kyqx,
\end{eqnarray} 
where $x$, $y\in (A,q)$, and $k$ is a fixed scalar in the field $\mathbf{k}$.
\end{definition}

\medskip
The square bracket $[ \, , \, ] _{i}$ with $3\le i\le 6$ is also denoted by 
$[ \, , \, ] _{i, k}$ to indicate its dependence on the scalar $k$.
The next proposition gives the basic property of the six square brackets.

\medskip
Each $i$-th square bracket with $6\ge i\ge 1$ can be expressed in terms of the Hu-Liu product $\circ_i$ by a few ways. One of the ways is given by
$$
[x, y ]_1=x\circ_{_{2,0}}y-y\circ_{_{2,0}}x,\qquad\quad [x, y ]_2=x\circ_{_{8,0,1}}y-y\circ_{_{8,0,1}}x,
$$
$$
[x, y ]_{3, k}=x\circ_{_{8,1,k}}y-y\circ_{_{8,1,k}}x,\qquad\quad 
[x, y ]_{4, k}=x\circ_{_{10,k}}y-y\circ_{_{10,k}}x,
$$
$$
[x, y ]_{5, k}=x\circ_{_{13,k}}y-y\circ_{_{13,k}}x,\qquad\quad 
[x, y ]_{6, k}=x\circ_{_{14,k}}y-y\circ_{_{14,k}}x.
$$

\medskip
\begin{proposition}\label{pr01.4} Let $(A,q)$ be an invariant algebra.
\begin{description}
\item[(i)] The $1$-st and $2$-nd square brackets  satisfy the Jacobi identity; that is, 
\begin{equation}\label{eq01.80}
[[x,y]_i,z]_i+[[y,z]_i,x]_i+[[z,x]_i,y]_i=0 \qquad\mbox{for $x,y,z\in (A,q)$,}
\end{equation}
where $i=1$ and $2$.
\item[(ii)] The $3$-rd square bracket  satisfies the following {\bf long Jacobi-like identity$^{1-st}$}:
\begin{eqnarray}\label{eq01.080}
&&[[x, y]_h, z]_k+[[y, z]_h, x]_k+[[z, x]_h, y]_k+\nonumber\\
&+& [[x, y]_k, z]_h+[[y, z]_k, x]_h+[[z, x]_k, y]_h=0
\end{eqnarray}
for $x,y,z\in (A,q)$ and $h$, $k\in \mathbf{k}$.
\item[(iii)]  Let $x,y,z\in (A,q)$ and $k$, $h\in \mathbf{k}$. If $i=4$, $5$ and $6$, then the $i$-th angle bracket  satisfies the 
{\bf Jacobi-like identity$^{1-st}$}:
\begin{equation}\label{eq01.801}
[[ x,  y] _{i, h}, z] _{i, k}
+[[ y, z] _{i, h}, x ] _{i, k}+ [[ z,  x] _{i, h}, y] _{i, k}=0 .
\end{equation}
Moreover, we have
\begin{equation}\label{eq01.802}
[[ x,  y] _{i, h}, z] _{i, k}=[[ x,  y] _{i, k}, z] _{i, h}\qquad\mbox{for $i=4$, $5$ and $6$.}
\end{equation}
\end{description}
\end{proposition}

\medskip
Since the six square brackets are anti-commutative, each of the six square brackets makes an invariant algebra into a Lie algebra.

\medskip
\section{Jordan Algebras}

\medskip
We begin this section by recalling the definition of a Jordan algebra from \cite{N}.

\begin{definition}\label{def6.4} A {\bf Jordan algebra} $J$ is an algebra over a field $\mathbf{k}$ of characteristic $\ne 2$ with a product composition $\odot$  satisfying 
\begin{equation}\label{eq6.263} 
x\odot y=y\odot x \qquad\mbox{(Commutative Law)}
\end{equation}
and
\begin{equation}\label{eq6.264} 
\Big((x\odot x)\odot y\Big) \odot x= (x\odot x)\odot (y \odot x)\qquad\mbox{(Jordan Identity)}
\end{equation}
for all $x$, $y\in J$.
\end{definition}

Based on the  $i$-th square brackets in Definition 3.1, we have the following

\begin{proposition}\label{pr6.4} An invariant algebra $(A, q)$ over a field $\mathbf{k}$ of characteristic $\ne 2$ is a Jordan algebra under each of the following $6$ products:
\begin{eqnarray}
\label{eq6.265} x\odot _{_1}y &=&\frac12 (qxy+qyx),\\
\label{eq6.266} x\odot _{_2}y&=&\frac12  (xqy+yqx),\\
\label{eq6.267} x\odot _{_3}y&=&\frac12  (xy+kxqy+yx+kyqx),\\
\label{eq6.268} x\odot _{_4}y&=&\frac12  (xy-xqy-kqxy+yx-yqx-kqyx),\\
\label{eq6.269} x\odot _{_5}y&=&\frac12  (xy-xyq-kqxy+yx-yxq-kqyx),\\
\label{eq6.270} x\odot _{_6}y&=&\frac12  (xy-xyq+kxqy+yx-yxq+kyqx), 
\end{eqnarray}
where $x$, $y\in (A,q)$, and $k$ is a fixed scalar in the field $\mathbf{k}$.
\end{proposition}

\medskip
\section{Leibniz Algebras}

\medskip
We begin this section with the definition of $i$-th  angle brackets.

\begin{definition}\label{def01.3} Let $(A,q)$ be an invariant algebra over a field $\mathbf{k}$. For $i=1$, $2$, $3$ and $4$, the {\bf $i$-th angle bracket} $\langle , \rangle _i$ is defined by
\begin{eqnarray}
\label{eq01.56}\langle x, y \rangle _1 :&=& xqy-qyx,\\
\label{eq01.57}\langle x, y \rangle _{2} :&=&xy-yx+yqx-xyq+k qyx-k qxy,\\
\label{eq01.58}\langle x, y \rangle _{3} :&=&xy-yx-xyq+yxq +k xqy -k qyx,\\
\label{eq01.59}\langle x, y \rangle _{4} :&=&xy-yx+yqx-xyq+k xqy -k qyx,
\end{eqnarray} 
where $x$, $y\in (A,q)$, and $k$ is a fixed scalar in the field $\mathbf{k}$.
\end{definition}

The angle bracket $\langle \, , \, \rangle _{i}$ with $2\le i\le 4$ is also denoted by $\langle \, , \, \rangle _{i, k}$ to indicate its dependence on the scalar $k$. The next proposition gives the basic property of the four angle brackets.

\begin{proposition}\label{pr01.3} Let $(A,q)$ be an invariant algebra.
\begin{description}
\item[(i)] The $1$-st angle bracket $\langle x, y \rangle _1$ satisfies the (right) Leibniz identity; that is, 
\begin{equation}\label{eq01.60}
\langle x, \langle y, z\rangle _1 \rangle _1=\langle\langle x,  y\rangle _1, z\rangle _1 - \langle\langle x,  z\rangle _1, y\rangle _1
\end{equation}
for $x,y,z\in (A,q)$.
\item[(ii)]  Let $x,y,z\in (A,q)$ and $k$, $h\in \mathbf{k}$. If $i=2$, $3$ and $4$, then the $i$-th angle bracket  satisfies the {\bf Jacobi-like identity$^{2-nd}$}:
\begin{equation}\label{eq01.601}
\langle x, \langle y, z\rangle _{i, k} \rangle _{i, h}=
\langle\langle x,  y\rangle _{i, h}, z\rangle _{i, k}
 - \langle\langle x,  z\rangle _{i, k}, y\rangle _{i, h} .
\end{equation}
Moreover, we have
\begin{equation}\label{eq01.602}
\langle x, \langle y, z\rangle _{i, k} \rangle _{i, h}=
\langle x, \langle y, z\rangle _{i, h} \rangle _{i, k}\qquad\mbox{for $i=2$, $3$ and $4$}
\end{equation}
and
\begin{equation}\label{eq01.603}
\langle\langle x,  y\rangle _{i, k}, z\rangle _{i, h}=
\langle\langle x,  y\rangle _{i, h}, z\rangle _{i, k}\qquad\mbox{for $i=2$ and $3$.}
\end{equation}
\end{description}
\end{proposition}

\medskip
Following \cite{Loday}, the notion of (right) Leibniz algebras is given in the following

\medskip
\begin{definition}\label{def1.1} A vector space $L$ is called a {\bf (right) Leibniz algebra} if there exists a binary operation $\langle \, , \, \rangle$: $L\times L\to L$, called the {\bf angle bracket},  such that the {\bf (right) Leibniz identity} holds:
\begin{equation}\label{eq1}
\langle\langle x, y\rangle , z\rangle=\langle x, \langle y, z\rangle\rangle+ 
\langle\langle x,  z\rangle , y\rangle \qquad\mbox{for $x,y,z\in L$.}
\end{equation}
\end{definition}

\medskip
By Proposition 5.1, each of the four angle brackets makes an invariant algebra into a (right) Leibniz algebra. 

\medskip
\section{Pre-Lie Algebras}

We begin this section by recalling the definition of a pre-Lie algebra from \cite{B}.

\begin{definition}\label{def6.2} A {\bf pre-Lie algebra} $A$ is a vector space equipped with a {\bf dot operation} $\cdot\, : A\otimes A\to A$ which satisfy the following identity:
\begin{equation}\label{eq6.175} 
(x\cdot y)\cdot  z-x\cdot (y\cdot  z)=(x\cdot z)\cdot  y-x\cdot (z\cdot  y)\quad\mbox{for $x$, $y$, $z\in A$.}
\end{equation}
\end{definition}

\medskip
We use $(A, \cdot)$ to denote a pre-Lie algebra $A$ equipped with a dot operation $\cdot$. Clearly, the square bracket
\begin{equation}\label{eq6.176}   
(x, y)\mapsto [x,\, y]^{\cdot}:=x\cdot y -y\cdot x \quad\mbox{for $x$, $y\in A$}
\end{equation}
satisfies the Jacobi identity; that is, a pre-Lie algebra is a Lie-admissible algebra. The square bracket $[x,\, y]^{\cdot}$ defined by (\ref{eq6.176}) is called the {\bf accompanying square product}. The following proposition gives $14$ ways of introducing a pre-Lie algebra structure in an invariant algebra and shows that each corresponding accompanying square product only differs from a square product in Definition 3.1 by a scalar.

\begin{proposition}\label{pr6.2} Let $(A, q)$ be an invariant algebra over a field $\mathbf{k}$. Let $k$ be a fixed scalar in the field $\mathbf{k}$. $((A, q), \cdot_{_i})$ is a pre-Lie algebra if the dot operation 
$\cdot_{_i}$ is chosen in one of the following $14$ ways:
\begin{description}
\item[(1)] $x\cdot_{_1}y=kyx+xqy-kyqx$, $[x,\, y]^{\cdot_{_1}}=[x,\, y]_{{_2}}$ for $k=0$ and 
$-\displaystyle\frac{1}{k}[x,\, y]^{\cdot_{_1}}=[x,\, y]_{{_3}}$ for $k\ne 0$;
\item[(2)] $x\cdot_{_2}y=kyx+xqy-kyqx-qyx-qxy$, $[x,\, y]^{\cdot_{_2}}=[x,\, y]_{{_2}}$ for $k=0$ and $-\displaystyle\frac{1}{k}[x,\, y]^{\cdot_{_2}}=[x,\, y]_{{_3}}$ for $k\ne 0$;
\item[(3)] $x\cdot_{_3}y=kyx+xqy+(1-k)yqx-qyx$, $[x,\, y]^{\cdot_{_3}}=[x,\, y]_{{_1}}$ for $k=0$ and $-\displaystyle\frac{1}{k}[x,\, y]^{\cdot_{_3}}=[x,\, y]_{{_4}}$ for $k\ne 0$;
\item[(4)] $x\cdot_{_4}y=kyx+xqy-qyx-qxy-kyxq$, $[x,\, y]^{\cdot_{_4}}=[x,\, y]_{{_2}}$ for $k=0$ and $-\displaystyle\frac{1}{k}[x,\, y]^{\cdot_{_4}}=[x,\, y]_{{_6}}$ for $k\ne 0$;
\item[(5)] $x\cdot_{_5}y=kyx+xqy+yqx-qyx-kyxq$, $-[x,\, y]^{\cdot_{_5}}=[x,\, y]_{{_1}}$ for $k=0$ and $-\displaystyle\frac{1}{k}[x,\, y]^{\cdot_{_5}}=[x,\, y]_{{_5}}$ for $k\ne 0$;
\item[(6)] $x\cdot_{_6}y=yx-xqy-yqx+xyq$ and $-[x,\, y]^{\cdot_{_6}}=[x,\, y]_{{_5}}$;
\item[(7)] $x\cdot_{_7}y=yx+kxqy-yqx+xyq$ and $-[x,\, y]^{\cdot_{_7}}=[x,\, y]_{{_6}}$;
\item[(8)] $x\cdot_{_8}y=yx+kxqy-yqx+xyq-(k+1)qyx-(k+1)qxy$ and 
$-[x,\, y]^{\cdot_{_8}}=[x,\, y]_{{_6}}$;
\item[(9)] $x\cdot_{_9}y=xy+kyqx-xyq-yxq$ and $[x,\, y]^{\cdot_{_9}}=[x,\, y]_{{_3}}$;
\item[(10)] $x\cdot_{_{10}}y=xy+yqx-xyq+kqyx-yxq$ and 
$[x,\, y]^{\cdot_{_{10}}}=[x,\, y]_{{_4}}$;
\item[(11)] $x\cdot_{_{11}}y=xy+yqx-xyq+kqxy-yxq$ and 
$[x,\, y]^{\cdot_{_{11}}}=[x,\, y]_{{_4}}$;
\item[(12)] $x\cdot_{_{12}}y=xy+kxqy-xyq-kqyx-kqxy$ and 
$[x,\, y]^{\cdot_{_{12}}}=[x,\, y]_{{_6}}$;
\item[(13)] $x\cdot_{_{13}}y=xy+kxqy+kyqx-xyq-kqyx$ and 
$[x,\, y]^{\cdot_{_{13}}}=[x,\, y]_{{_5}}$;
\item[(14)] $x\cdot_{_{14}}y=xy+kxqy-(k+1)qyx-(k+1)qxy$ and 
$[x,\, y]^{\cdot_{_{14}}}=[x,\, y]_{{_3}}$,
\end{description} 
where $x$, $y\in (A, q)$, and $[x, \, y]_{_i}$ is the $i$-th square product.
\end{proposition}

\medskip
\section{Left Symmetric Algebras}

\medskip
We begin this section by recalling the definition of a left symmetric algebra from \cite{B}.

\medskip
\begin{definition}\label{def6.3} A {\bf left symmetric algebra} $A$ is a vector space equipped with a {\bf dot operation} $\cdot\, : A\otimes A\to A$ which satisfy the following identity:
\begin{equation}\label{eq6.219} 
x\cdot (y\cdot  z)-(x\cdot y)\cdot  z=y\cdot (x\cdot  z)-(y\cdot x)\cdot  z\quad\mbox{for $x$, $y$, $z\in A$.}
\end{equation}
\end{definition}

We use $(A, \cdot)$ to denote a left symmetric algebra $A$ equipped with a dot operation $\cdot$. It is easy to check that the square bracket
\begin{equation}\label{eq6.220}   
(x, y)\mapsto [x,\, y]^{\cdot}:=x\cdot y -y\cdot x \quad\mbox{for $x$, $y\in A$}
\end{equation}
satisfies the Jacobi identity; that is, a left symmetric algebra is a Lie-admissible algebra. The square bracket $[x,\, y]^{\cdot}$ defined by (\ref{eq6.220}) is called the {\bf accompanying square product}. The following proposition gives $14$ ways of introducing a left symmetric algebra structure in an invariant algebra and shows that each corresponding accompanying square product only differs from a  square product in Definition 3.1 by a scalar.

\begin{proposition}\label{pr6.3} Let $(A, q)$ be an invariant algebra over a field $\mathbf{k}$. Let $k$ be a fixed scalar in the field $\mathbf{k}$. $((A, q), \cdot_{_i})$ is a left symmetric algebra if the dot operation 
$\cdot_{_i}$ is chosen in one of the following $14$ ways:
\begin{description}
\item[(1)] $x\cdot_{_1}y=kxy-kxqy+yqx$, $-[x,\, y]^{\cdot_{_1}}=[x,\, y]_{{_2}}$ for $k=0$ and 
$\displaystyle\frac{1}{k}[x,\, y]^{\cdot_{_1}}=[x,\, y]_{{_3}}$ for $k\ne 0$;
\item[(2)] $x\cdot_{_2}y=kxy+(1-k)xqy+yqx-qxy$, $-[x,\, y]^{\cdot_{_2}}=[x,\, y]_{{_1}}$ for $k=0$ and 
$\displaystyle\frac{1}{k}[x,\, y]^{\cdot_{_2}}=[x,\, y]_{{_4}}$ for $k\ne 0$;
\item[(3)] $x\cdot_{_3}y=xy+yqx-qxy$ and $[x,\, y]^{\cdot_{_3}}=[x,\, y]_{{_4}}$;
\item[(4)] $x\cdot_{_4}y=kxy-kxqy+yqx-qyx-qxy$, $-[x,\, y]^{\cdot_{_4}}=[x,\, y]_{{_2}}$ for $k=0$ and 
$\displaystyle\frac{1}{k}[x,\, y]^{\cdot_{_4}}=[x,\, y]_{{_3}}$ for $k\ne 0$;
\item[(5)] $x\cdot_{_5}y=kyx+(1-k)yqx-qyx-qxy$, $-[x,\, y]^{\cdot_{_5}}=[x,\, y]_{{_2}}$ for $k=0$ and 
$\displaystyle\frac{1}{k}[x,\, y]^{\cdot_{_5}}=[x,\, y]_{{_3}}$ for $k\ne 0$;
\item[(6)] $x\cdot_{_6}y=kyx+yqx-qyx-qxy-kyxq$, $-[x,\, y]^{\cdot_{_6}}=[x,\, y]_{{_2}}$ for $k=0$ and 
$-\displaystyle\frac{1}{k}[x,\, y]^{\cdot_{_6}}=[x,\, y]_{{_6}}$ for $k\ne 0$;
\item[(7)] $x\cdot_{_7}y=kyx+xqy+yqx-qxy-kyxq$, $-[x,\, y]^{\cdot_{_7}}=[x,\, y]_{{_1}}$ for $k=0$ and 
$-\displaystyle\frac{1}{k}[x,\, y]^{\cdot_{_7}}=[x,\, y]_{{_5}}$ for $k\ne 0$;
\item[(8)] $x\cdot_{_8}y=kxy+yqx-kxyq-qyx-qxy$, $-[x,\, y]^{\cdot_{_8}}=[x,\, y]_{{_2}}$ for $k=0$ and 
$\displaystyle\frac{1}{k}[x,\, y]^{\cdot_{_8}}=[x,\, y]_{{_6}}$ for $k\ne 0$;
\item[(9)] $x\cdot_{_9}y=kxy+xqy+yqx-kxyq-qxy$, $-[x,\, y]^{\cdot_{_9}}=[x,\, y]_{{_1}}$ for $k=0$ and 
$\displaystyle\frac{1}{k}[x,\, y]^{\cdot_{_9}}=[x,\, y]_{{_5}}$ for $k\ne 0$;
\item[(10)] $x\cdot_{_{10}}y=yx+xqy-xyq+kqxy-yxq$ and 
$[x,\, y]^{\cdot_{_{10}}}=[x,\, y]_{{_4}}$;
\item[(11)] $x\cdot_{_{11}}y=yx+kxqy-xyq-yxq$ and 
$-[x,\, y]^{\cdot_{_{11}}}=[x,\, y]_{{_3}}$;
\item[(12)] $x\cdot_{_{12}}y=yx+kxqy+(k-1)yqx-xyq-(k-1)qxy-yxq$ and 
$-[x,\, y]^{\cdot_{_{12}}}=[x,\, y]_{{_4}}$;
\item[(13)] $x\cdot_{_{13}}y=yx+xqy-xyq+kqyx-yxq$ and 
$-[x,\, y]^{\cdot_{_{13}}}=[x,\, y]_{{_4}}$;
\item[(14)] $x\cdot_{_{14}}y=yx+xqy+kyqx-xyq-kqyx-kqxy-yxq$ and 
$-[x,\, y]^{\cdot_{_{14}}}=[x,\, y]_{{_3}}$,
\end{description} 
where $x$, $y\in (A, q)$, and $[x, \, y]_{_i}$ is the $i$-th square product.
\end{proposition}

\medskip
\section{Representation$^{(c_1, \dots, c_8)}$ of Algebras }

Let $(\mathcal{A}, \star)$ be an algebra (not necessarily associative algebra) over  $\mathbf{k}$ with a product $\star$. Based on invariant algebras, we introduce  a representation$^{(c_1, \dots, c_8)}$ of an algebra  $(\mathcal{A}, \star)$ in the following

\begin{definition}\label{def7.1} Let $W$ be a subspace of a vector space $V$ over a field $\mathbf{k}$ and let $q$ be a $W$-idempotent. A linear map $\varphi$ from an algebra  
$(\mathcal{A}, \star)$ to the invariant algebra $(End(V), q)$ is called a {\bf representation$^{(c_1, \dots, c_8)}$} of $(\mathcal{A}, \star)$ induced by $(q, W)$ if there exist the scalars $c_1, \dots, c_8\in \mathbf{k}$ such that
\begin{eqnarray}\label{eq7.1}
&&\varphi(x\star y)
=c_1\varphi(x)\varphi(y)+c_2\varphi(y)\varphi(x)+c_3q\varphi(x)\varphi(y)
+c_4q\varphi(y)\varphi(x)+\nonumber\\
&&\quad +c_5\varphi(x)q\varphi(y)+c_6\varphi(y)q\varphi(x)
+c_7\varphi(x)\varphi(y)q+c_8\varphi(y)\varphi(x)q
\end{eqnarray}
for all $x$, $y\in \mathcal{A}$.
\end{definition}

The language of  modules$^{(c_1, \dots, c_8)}$ is more convenient to state some facts about representations$^{(c_1, \dots, c_8)}$. We now introduce a module$^{(c_1, \dots, c_8)}$ over an algebra  in the following

\begin{definition}\label{def7.2} Let $W$ be a subspace of a vector space $V$ over a field $\mathbf{k}$ and let $q$ be a $W$-idempotent. $V$ is called a 
{\bf module$^{(c_1, \dots, c_8)}$} over an algebra $(\mathcal{A}, \star)$ or a {\bf $\mathcal{A}$-module$^{(c_1, \dots, c_8)}$} induced by $(q, W)$ if there is a map: 
$(x, v)\mapsto x\cdot v$ from $\mathcal{A}\times V$ to $V$ such that
\begin{equation}\label{eq7.2} 
(ax+by)\cdot v=a(x\cdot v)+b(y\cdot v),
\end{equation}
\begin{equation}\label{eq7.3} 
x\cdot (av+bu)=a(x\cdot v)+b(x\cdot u),
\end{equation}
\begin{eqnarray}\label{eq7.4}
&&(x\star y)\cdot v
=c_1x\cdot y\cdot v+c_2y\cdot x\cdot v+c_3q\cdot x\cdot y\cdot v
+c_4q\cdot y\cdot x\cdot v+\nonumber\\
&&\quad +c_5x\cdot q\cdot y\cdot v+c_6y\cdot q\cdot x\cdot v
+c_7x\cdot y\cdot q\cdot v+c_8y\cdot x\cdot q\cdot v,
\end{eqnarray}
\begin{equation}\label{eq7.5} 
x\cdot W:=\{\, x\cdot w \,|\, w\in W\,\}\subseteq W,
\end{equation}
where $x$, $y\in \mathcal{A}$, $v$, $u\in V$, $a$, $b\in\mathbf{k}$ and $q\cdot v:=q(v)$.
\end{definition}

A  module$^{(c_1, \dots, c_8)}$ over an algebra $(\mathcal{A}, \star)$ induced by $(q, W)$ is also denoted by $V_{(q, W)}$. A subspace $U$ of an $\mathcal{A}$-module$^{(c_1, \dots, c_8)}$
$V=V_{(q, W)}$ is called a {\bf submodule$^{(c_1, \dots, c_8)}$} of $V$ if
\begin{equation}\label{eq7.6} 
x\cdot u\in U \quad\mbox{for $x\in\mathcal{A}$ and $u\in U$}
\end{equation}
and
\begin{equation}\label{eq7.7} 
q\cdot u-u\in U \quad\mbox{for $u\in U$.}
\end{equation}

Every $\mathcal{A}$-module$^{(c_1, \dots, c_8)}$ $V=V_{(q, W)}$ has at least three  submodules$^{(c_1, \dots, c_8)}$: $0$, $W$ and $V$. An 
$\mathcal{A}$-module$^{(c_1, \dots, c_8)}$ with $c_3=\dots =c_8=0$ is called an 
{\bf $\mathcal{A}$-module$^{{c_1 \choose c_2}}$}. 

\begin{proposition}\label{pr7.1} If $V_{(q, W)}$ is module$^{(c_1, \dots, c_8)}$ over an algebra $(\mathcal{A}, \star)$, then the subspace $W$ is an $\mathcal{A}$-module$^{{c_1 \choose c_2}}$ by restriction and the quotient space $\displaystyle\frac{V}{W}$ is a 
$\mathcal{A}$-module$^{{c_1+c_3+c_5+c_7 \choose c_2+c_4+c_6+c_8}}$ under the following action:
\begin{equation}\label{eq7.8} 
x\cdot (v+W):=x\cdot v +W \quad\mbox{for $x\in\mathcal{A}$ and $v\in V$.}
\end{equation}
\end{proposition}

\bigskip
We now define the building blocks for modules$^{(c_1, \dots, c_8)}$ over an algebra.

\begin{definition}\label{def7.3} Let $V_{(q, W)}$ be a module$^{(c_1, \dots, c_8)}$ over an algebra $(\mathcal{A}, \star)$.
\begin{description}
\item[(i)] $V_{(q, W)}$ is said to be {\bf $2$-irreducible} if $V_{(q, W)}\ne 0$ and 
$V_{(q, W)}$ has no submodules$^{(c_1, \dots, c_8)}$ which are not equal to $0$ and 
$V_{(q, W)}$.
\item[(ii)] $V_{(q, W)}$ is said to be {\bf $3$-irreducible} if $W\ne 0$, $V_{(q, W)}\ne W$ and $V_{(q, W)}$ has no submodules$^{(c_1, \dots, c_8)}$ which are not equal to $0$, $W$ and 
$V_{(q, W)}$.
\end{description}
\end{definition}

The next proposition gives the basic properties of the building blocks.

\begin{proposition}\label{pr7.2} Let $V=V_{(q, W)}$ be a module$^{(c_1, \dots, c_8)}$ over an algebra $(\mathcal{A}, \star)$.
\begin{description}
\item[(i)] If $V_{(q, W)}$ is  $2$-irreducible, then either $V=W$ in which case $q=0$ on $V$ and $V$ is a  $\mathcal{A}$-module$^{{c_1 \choose c_2}}$, or $W=0$ in which case $q=1$ on $V$ and $V$ is a $\mathcal{A}$-module$^{{c_1+c_3+c_5+c_7 \choose c_2+c_4+c_6+c_8}}$.
\item[(ii)] If $V_{(q, W)}$ is $3$-irreducible, then both the 
$\mathcal{A}$-module$^{{c_1 \choose c_2}}$ $W$ and the 
$\mathcal{A}$-module$^{{c_1+c_3+c_5+c_7 \choose c_2+c_4+c_6+c_8}}$  $\displaystyle\frac{V}{W}$ are $2$-irreducible.
\end{description}
\end{proposition}

\begin{definition}\label{def7.2} Let $(V_1)_{(q_1, W_1)}$ and $(V_2)_{(q_2, W_2)}$ be two modules$^{(c_1, \dots, c_8)}$ over an algebra $(\mathcal{A}, \star)$. 
\begin{description}
\item[(i)] A linear map $\varphi : V_1\to V_2$ is called a 
{\bf homomorphism$^{(c_1, \dots, c_8)}$} if
\begin{equation}\label{eq7.13} 
\varphi (x\cdot v_1)=x\cdot \varphi (v_1) \quad\mbox{for $x\in L$ and $v_1\in V_1$}
\end{equation}
and
\begin{equation}\label{eq7.14} 
\varphi (q_1\cdot v_1)=q_2\cdot \varphi (v_1) \quad\mbox{for $v_1\in V_1$.}
\end{equation}
\item[(ii)] If a homomorphism$^{(c_1, \dots, c_8)}$ $\varphi : V_1\to V_2$ is bijective, then we say that $\varphi$ is an {\bf isomorphism$^{(c_1, \dots, c_8)}$} and $V_1$ is {\bf isomorphic} to $V_2$.
\end{description}
\end{definition}

Let $V=V_{(q, W)}$ be a module$^{(c_1, \dots, c_8)}$ over an algebra $(\mathcal{A}, \star)$ induced by $(q, W)$. If $U$ is a  submodule$^{(c_1, \dots, c_8)}$ of $V$, then $q|U$ is a $(W\cap U)$-idempotent and $\bar{q}: \displaystyle\frac{V}{U}\to \displaystyle\frac{V}{U}$ is a $\displaystyle\frac{U+W}{U}$-idempotent, where $\bar{q}$ is defined by
\begin{equation}\label{eq7.15} 
\bar{q}\cdot (v+U)=\bar{q}\cdot v +U \quad\mbox{for $v\in V$.}
\end{equation}
Hence, $U_{q|U, W\cap U}$ is an $\mathcal{A}$-module$^{(c_1, \dots, c_8)}$ induced by 
$(q|U, W\cap U)$, and the quotient space  $\displaystyle\frac{V}{U}$ is also a $\mathcal{A}$-module$^{(c_1, \dots, c_8)}$ under the module$^{(c_1, \dots, c_8)}$ action consisting of (\ref{eq7.15}) and the following
\begin{equation}\label{eq7.16} 
x\cdot (v+U):=x\cdot v +U \quad\mbox{for $x\in\mathbf{A}$ and $v\in V$.}
\end{equation}
The $\mathcal{A}$-module$^{(c_1, \dots, c_8)}$ 
$\left(\displaystyle\frac{V}{U}\right)_{\left(\bar{q}, \frac{U+W}{U}\right)}$ is called the 
{\bf quotient $\mathcal{A}$-module$^{(c_1, \dots, c_8)}$} with respect to 
submodule$^{(c_1, \dots, c_8)}$ $U$.

\begin{proposition}\label{pr7.3} Let $(V_1)_{(q_1, W_1)}$ and $(V_2)_{(q_2, W_2)}$ be two modules$^{(c_1, \dots, c_8)}$ over an algebra $(\mathcal{A}, \star)$, and suppose that 
$\varphi : V_1\to V_2$ is a homomorphism$^{(c_1, \dots, c_8)}$.
\begin{description}
\item[(i)] $\varphi (W_1)\subseteq  W_2$.
\item[(ii)] The {\bf kernel} $Ker\varphi: =\{\, v_1\in V_1\,|\, \varphi (v_1)=0 \,\}$ is a 
submodule$^{(c_1, \dots, c_8)}$ of $\mathcal{A}$-module$^{(c_1, \dots, c_8)}$ $V_1$.
\item[(iii)] The {\bf image} $Im\varphi: =\{\, \varphi (v_1) \,|\, v_1\in V_1 \,\}$ is a 
submodule$^{(c_1, \dots, c_8)}$ of $\mathcal{A}$-module$^{(c_1, \dots, c_8)}$ $V_2$.
\item[(iv)] The map $\bar{\varphi}: \displaystyle\frac{V_1}{Ker\varphi}\to Im\varphi$ defined by
\begin{equation}\label{eq7.17} 
\varphi (v_1+Ker\varphi):=\varphi (v_1) \quad\mbox{for $v_1\in V_1$}
\end{equation}
is an isomorphism$^{(c_1, \dots, c_8)}$ from the quotient 
$\mathcal{A}$-module$^{(c_1, \dots, c_8)}$
$\left(\displaystyle\frac{V_1}{Ker\varphi}\right)$ induced by 
$\left(\bar{q},\, \frac{Ker\varphi +W_1}{Ker\varphi}\right)$ to the 
submodule$^{(c_1, \dots, c_8)}$ $U_{(q|U, W\cap U)}$.
\end{description}
\end{proposition}

\end{document}